\documentclass[11pt]{article}

\usepackage{latexsym,mathrsfs}
\usepackage{amsmath,amssymb} 
\usepackage{amsthm,enumerate,verbatim}
\usepackage{amsfonts}
\usepackage{graphicx}
\usepackage{algorithm}
\usepackage{algorithmic}
\usepackage{url}
\usepackage{color}
\usepackage{hyperref}
\usepackage{authblk}

\newcommand {\mat}  [1] {\left[\begin{array}{#1}}
\newcommand {\rix}      {\end{array}\right]}

\newtheorem{theorem}{Theorem}

\newtheorem{definition}{Definition}
\newtheorem*{definition*}{Definition}

\newtheorem{remark}{Remark}
\newtheorem*{problem*}{Problem}

\DeclareMathOperator{\sym}{sym} 
\DeclareMathOperator{\skewp}{skew}

\def \R{{\mathbb R}}
\def \C{{\mathbb C}}
\newcommand{\Tc}{\mathcal{T}}

\definecolor{brightpink}{rgb}{1.0, 0.0, 0.5}

\newcommand{\revise}[1]{{{\color{black} #1}}}

\newcommand{\reviseELA}[1]{{{\color{black} #1}}}

\title{Computing the nearest scattering passive system} 
\date{}

\author{
Karim Cherifi\thanks{Département AS2M, 
Équipe DATA-PHM, SUPMICROTECH-ENSMM,  
26 rue de l'Epitaphe, 
25000 BESANCON, Email: karim.cherifi@femto-st.fr} 
\qquad 
Nicolas Gillis\thanks{University of Mons, Rue de Houdain 9, 7000 Mons, Belgium. Email: nicolas.gillis@umons.ac.be.} \qquad 
 Punit Sharma\thanks{Indian Institute of Technology Delhi, Hauz Khas, 110016 New Delhi, India. Email: punit.sharma@maths.iitd.ac.in.}
 }

\begin{document}

\maketitle

\begin{abstract}
In this paper, we consider linear time-invariant control systems which are bounded real, also known as scattering passive. Our main theoretical contribution is to show the equivalence between such systems and port-Hamiltonian (PH) systems whose factors satisfy certain linear matrix inequalities. Based on this result, we propose a formulation for the problem of finding the nearest bounded real system to a given system, and design an algorithm combining alternating optimization and Nesterov's fast gradient method. This formulation also allows us to check whether a given system is bounded real by solving a semidefinite program, and provide a PH parametrization for it. We illustrate our proposed algorithms on real-world and synthetic data sets.  
\end{abstract}


\section{Introduction}

In this paper, we consider $m$-input $m$-output linear
time-invariant (LTI) control systems of the form
\begin{align}
\label{eq1:sys}
\begin{split}
\dot{x}(t) = & \; Ax(t)+Bu(t),
\\
y(t) = & \; Cx(t)+Du(t),
\end{split}
\end{align}
on the unbounded interval $t \in [t_0,\infty)$. Here, $A \in \R^{n,n}$, $B \in \R^{n,m}$, $C \in \R^{m,n}$, and $D \in \R^{m,m}$
are given matrices, $x(t)$ is the vector of state variables, $u(t)$ is the vector of inputs, and $y(t)$
is the vector of outputs. We use the matrix quadruple $(A,B,C,D)$ to refer to system in the form~\eqref{eq1:sys}.  The corresponding transfer function in frequency domain using the Laplace transform is written as 
\begin{equation}\label{tf}
\Tc(s)=C(sI-A)^{-1}B+D.
\end{equation}

%Port-Hamiltonian (PH) systems \cite{SchJ14,MehU23} have been introduced as generalization of Hamiltonian systems to open systems. Their structure allows for the consistent modeling of complex systems by ensuring that the interconnections are power-conserving making the PH framework particularly suitable for multiphysics modeling. In addition, PH systems are inherently passive and stable. 

%More particularly, we are interested in passive systems. 
The system~\eqref{eq1:sys} is said to be dissipative~\cite{Wil72} if there exists a  storage function $V: \R^{n} \rightarrow [0,\infty)$ such that $V(0)=0$ and the following dissipation inequality holds
	\begin{align}
		%\label{dissi_ineq}
	V(x(t_1)) - V(x(t_0)) \leq \int_{t_0}^{t_1} \revise{s}(y(\tau),u(\tau)) d \tau,  \nonumber
	\end{align}
	along all possible trajectories of the system starting at $x(t_0)$, for all $x(t_0)$, $t_1 \geq t_0$  and \revise{$s$}  is the supply rate. The quadratic supply rates of the following form are usually used for linear systems.
\begin{equation}
	\label{eq:supply}
	\mathcal{S}(y,u):=\begin{bmatrix}
		y\\ u \end{bmatrix}^\top \begin{bmatrix} N & \Omega\\\Omega^\top & M\end{bmatrix}\begin{bmatrix}
		y\\ u
	\end{bmatrix},\quad N,\Omega,M\in\mathbb{R}^{m,  m},\quad N=N^\top ,\quad M=M^\top . 
\end{equation}
In particular, we are interested in two special supply rates. The \emph{impedance supply rate} of the form 
\[ 
    \mathcal{S}_{imp}(y,u):= \begin{bmatrix}
		y\\ u \end{bmatrix}^\top \begin{bmatrix} 0 &  I_m\\ I_m & 0\end{bmatrix}\begin{bmatrix}
		y\\ u
	\end{bmatrix}= 2 y^\top  u, \nonumber
\] %\label{impsupply} 
and the \emph{scattering supply rate} 
\begin{equation}%\label{scatsupply}
	\mathcal{S}_{sca}(y,u):=\begin{bmatrix}
		y\\ u \end{bmatrix}^\top \begin{bmatrix} -I_m & 0\\0& I_m\end{bmatrix}\begin{bmatrix}
		y\\ u
	\end{bmatrix}=\|u\|^2-\|y\|^2, \nonumber
\end{equation}
where $I_m$ is the identity matrix of dimension $m$.  
Dissipative systems with an impedance supply rate are called \emph{impedance passive} and systems with a scattering supply rate are called \emph{scattering passive}.  

Passive systems are closely related to port-Hamiltonian (PH) systems, which have been introduced as a generalization of Hamiltonian systems to open systems; see~\cite{SchJ14,MehU23,RasFVS20} and references therein. Their structure allows for the consistent modeling of complex systems by ensuring that the interconnections are power-conserving making the PH framework particularly suitable for multiphysics modeling. A linear time-invariant input-state-output PH system can be written as 
\begin{align}
	\label{eq:phsystem1}
	\begin{split}
		\dot{x}(t) = & \; (J-R)Qx(t) + (F-P)u(t), \\
		y(t) = & \; (F+P)^\top   Qx(t) + Du(t),
	\end{split}
\end{align}
where $x$ is the state, $u$ the input, and $y$ the output, and the matrix $Q \in \R^{n,n}$ is symmetric positive definite. The function $x\rightarrow \frac{1}{2}x^\top   Q x$ is the \emph{Hamiltonian} and describes the energy of the system; the matrix $J^\top  =-J \in \R^{n,n}$ is the structure matrix that describes flux among energy storage elements; the matrix $R \in \R^{n,n}$ with $R$ being \revise{symmetric positive} semidefinite, is the dissipation matrix and describes the energy dissipation/loss in the system; and the matrices $F\pm P \in \R^{n,m}$ are the port matrices describing the manner in which energy enters and exists the system. A PH system in the form~\eqref{eq:phsystem1} is called impedance PH (iPH) if the matrices $R$, $P$ and $D$ satisfy
\begin{equation}\label{eq:posrealK}
K_i=\mat{cc}R &P\\P^\top   &\frac{D+D^\top}{2} \rix \succeq 0. 
\end{equation}
%i.e., $K_i$ is positive semidefinite.
In most of the literature on PH systems, the impedance supply rate is considered~\cite{CheGH23}. This supply rate has a physical meaning in many applications. For example the supply rate in electrical circuits represents power since the inputs and outputs are made up of the voltage and current. On the other hand, the scattering case is not broadly studied in the literature \cite{StaW12, CheGHM23,CheG23,CerVB07}. However, this supply rate is the most appropriate for many applications; for example in macro-modeling for electronic design~\cite{BraGZC22}. 
In that case, one is interested in the difference between the power put in the systems and power exiting the system which corresponds to the scattering supply rate. 
%Another example is the transport equation, e.g., \cite{TolMP24} where they study a boundary controlled transport equation which has non-collocated boundary control and observation, and reveals a scattering-energy preserving structure. In that case, the scattering supply rate represents the difference in energy between the input at the boundary and the output at another end of the boundary. Other applications include quantum mechanics and quantum field theory~\cite{LanL91, Bar67}.
\revise{
Another example is the transport equation, see, e.g., \cite{TolMP24}, where a boundary-controlled transport equation with non-collocated boundary control and observation has a scattering energy-preserving structure. In that case, the scattering supply rate represents the difference in energy between the input at one boundary and the output at the other. Other applications arise in quantum mechanics and quantum field theory~\cite{LanL91, Bar67}.} 

For minimal systems, impedance passivity is equivalent to being \emph{positive real} (PR), and scattering passivity  is equivalent to being \emph{bounded real} (BR)~\cite{BroLME20}. 
Motivated by the nearest PR system problem studied in~\cite{gillis2018finding}, in this paper, we aim to solve the nearest BR system problem defined as follows.
%. More precisely, 
%for a given $(A,B,C,D)$, the main goal of this paper is to solve the following problem.
\begin{problem*}[Nearest BR system] 
	For a given system $(A,B,C,D)$, find the \emph{nearest BR transfer function} with system matrices $(\tilde A,\tilde B,\tilde C, \tilde D)$ to $(A,B,C,D)$, that is, solve
	\begin{equation} \label{probdef}
		\inf_{(\tilde A,\tilde B,\tilde C, \tilde D) \in \mathcal{BR}}
		\mathcal{F}(\tilde A,\tilde B,\tilde C,\tilde D),
		\tag{$\mathcal P$}
	\end{equation}
	where
	\begin{equation} \label{def:FinP}
		\mathcal{F}(\tilde A,\tilde B,\tilde C,\tilde D) =
		{\|A-\tilde A\|}_F^2+{\|B-\tilde B\|}_F^2
		+{\|C-\tilde C\|}_F^2+{\|D-\tilde D\|}_F^2,
	\end{equation}
     ${\|\cdot\|}_F$ denotes the Frobenius norm, 
 and $\mathcal{BR}$ is the set of all BR systems
	$(\tilde A,\tilde B,\tilde C, \tilde D)$, where $\tilde A \in \R^{n,n}$, $\tilde B \in \R^{n,m}$, $\tilde C \in \R^{m,n}$, and $\tilde D \in \R^{m,m}$.
\end{problem*}
%
% The nearest BR system problem~\eqref{probdef} is motivated by the nearest PR system problem, studied in~\cite{gillis2018finding}, where the authors proposed an optimization based method to find a nearby PR system to a given non PR system by finding a nearby impedance passive PH system. 

Although the problem~\eqref{probdef} is in principle similar to the one defined in~\cite{gillis2018finding} for positive real transfer functions, the nearest BR problem leads to a different set of constraints and thus to a different optimization problem with its own set of challenges. A solution to the nearest BR system problem~\eqref{probdef} can be useful in passive PH representations of the system from data. 
Although PH representations are usually derived from first principles with a strong emphasis on the energy of the system, in many applications, only input/output data is available. System identification techniques are then used to recover the system matrices. System identification methods for PH systems include methods based on optimization~\cite{gillis2018finding,SchV23} and the Loewner framework~\cite{Che22}.
%, and positive real balanced truncation~\cite{CheMH19}. 
However, since most system identification methods do not enforce passivity, one needs to check if the resulting system is passive. If it is not, then a correction step is required to obtain the nearest passive model~\cite{gillis2018finding,FazGL21}; also called passivity enforcement \cite{Gri04,GriU07,GriG15,PorVIS08,Bjo21,RoI21} or passivation in the recent paper~\cite{nicodemus2025klap}. In particular, in this paper, we focus on finding the nearest scattering passive system.

\paragraph{Contribution and outline of the paper} 

The remainder of the paper is structured as follows. 
In Section~\ref{sec:pHBR}, we define scattering PH (sPH) systems, that is,  PH systems that are scattering passive. We then obtain a necessary and sufficient condition to show the equivalence between sPH and BR systems. For minimal sytems, this yields a new way to show the equivalence between scattering passivity and bounded realness (BRness). In Section~\ref{sec:NearsPH}, based on these conditions, we propose a formulation to solve the nearest BR system problem.
Our principle strategy for computing a nearby BR system to a given non BR system is to reformulate the highly non-convex problem into an equivalent (non-convex) problem with a convex feasible set onto which points can be projected relatively easily. 
 This leads us, in Section~\ref{sec:checkBR}, 
to propose a semidefinite programming (SDP) formulation to check whether a given system is BR. To handle large problems, we design an algorithm that combines alternating optimization and Nesterov's fast gradient method (Algorithm~\ref{alg:checkBRfgm}).    
In Section~\ref{sec:EAOnearestBRsPHform}, we use a similar algorithmic strategy to compute a nearby BR system in sPH form (Algorithm~\ref{alg:eao}). 
We illustrate the effectiveness of these algorithms on real-world and synthetic systems in Section~\ref{sec:numexp}. 

%\paragraph{Notation} 
%We denote by ${\|\cdot\|}_F$ the Frobenius norm and by $*$ the complex conjugate transpose. We write $A \succ 0$ (resp., $A \succeq 0$) if $A$ is symmetric symmetric symmetric positive definite (resp., semideﬁnite). 
%The real part of $s\in \C$  is denoted by $\operatorname { R e }(s)$ and $j$ stands for the imaginary number. 
%For a Hermitian matrix $H$, $\lambda_{\max}(H)$ and $\lambda_{\min}(H)$ denote the largest and the smallest eigenvalue of $H$, respectively. 

\section{PH systems and their relationship to BR systems} 
\label{sec:pHBR}

In this section, we study the link between BR systems, scattering passive systems and PH systems. 
\begin{definition}[\cite{BroLME20}]%\rm
A system $(A,B,C,D)$ in~\eqref{eq1:sys} is said to be BR  if its transfer function $\Tc(s)$ defined in~\eqref{tf}  satisfies the following two conditions:  
\begin{enumerate}
    \item[(a)]~all elements of $\Tc(s)$ are analytic for 
 $\operatorname { R e }(s) \geq 0$ where $\operatorname { R e }(s)$ is the real part of $s$, and 
 \item[(b)]~$I_m-(\Tc(j \omega))^* \Tc(j \omega) \succcurlyeq 0$ for all $\omega \in \mathbb{R}$, where 
 $^*$ denotes the complex conjugate transpose, and $j$ stands for the imaginary unit. 
\end{enumerate}
\end{definition}

By the Kalman–Yakubovich–Popov (KYP) or BR lemma for a system~\eqref{eq1:sys}, \revise{strict BRness (i.e., (b) with ``positive definite" instead of ``positive semidefinite")} can be characterized in terms of symmetric positive definite solutions $X$ of the following LMI~\cite{CavF19} : 
\begin{equation}\label{eq:LMIO}
\mat{cc} A^\top  X +X A +C^\top  C & XB+C^\top  D \\B^\top  X+D^\top  C & D^\top  D-I_m \rix \prec 0.
\end{equation}
Using the Schur complement, the LMI in~\eqref{eq:LMIO} can be equivalently written as $\revise{{BR}_{LMI}(X)} \prec 0$, where 
\begin{equation}\label{eq:LMI1strict}
{BR}_{LMI}(X):=\mat{ccc} A^\top  X +XA & XB & C^\top   \\B^\top  X & -I_m & D^\top   \\ C & D & -I_m \rix.
\end{equation}
For minimal systems $(A,B,C,D)$ in the form~\eqref{eq1:sys}, that is, when the system is controllable and observable, BRness is equivalent to scattering passivity, which 
%the LMI ~\eqref{eq:LMI1strict} can be relaxed by replacing $\prec$ with $\preceq$.
can be characterized in terms of symmetric positive definite solutions $X$ of the equation 
\[
{BR}_{LMI}(X) \preceq 0,
\]
where the matrix \revise{${BR}_{LMI}(X)$} is defined by~\eqref{eq:LMI1strict}~\cite{Ran96,CavF19}.

As mentioned earlier, the impedance PH systems have been studied extensively in the literature. On the other hand, the scattering PH systems have not received much attention. 
\begin{definition}
    A PH system in the form~\eqref{eq:phsystem1} is said to be \emph{scattering PH} (sPH) if its factor matrices $R,F,P$, and $D$ satisfy 
\begin{equation}\label{eq:bourealK}
K_s:=\mat{ccc}2R & -(F-P) & -(F+P)\\ -(F-P)^\top   & I_m & -D^\top  \\ -(F+P)^\top   & -D & I_m \rix \succeq 0.
\end{equation}
\end{definition}
Now we study the relation between sPH systems and BR starting with the case where $K_s \succ 0$. %\reviseELA{and prove two theorems leading to the equivalence of sPH  and BR systems.} 
\begin{theorem}\label{thm:sphtoscat}
    Every sPH system in the form~\eqref{eq:phsystem1} with $K_s \succ 0$ is \revise{strict BR}. 
\end{theorem}
\begin{proof} Consider a sPH system in the form~\eqref{eq:phsystem1} such that the matrices $R,F,P$, and $D$ satisfy that $K_s \succ 0$, where $K_s$ is as defined in~\eqref{eq:bourealK}. The result \reviseELA{follows by} using the fact that  $X=Q$ satisfies  ${BR}_{LMI}(X) \prec 0$, where ${BR}_{LMI}(X)$ is as defined in~\eqref{eq:LMI1strict} with $A=(J-R)Q$, $B=F-P$, $C=(F+P)^\top Q$. Indeed, we have
\begin{align*}
&\mat{ccc} ((J-R)Q)^\top  Q +Q^\top   (J-R)Q & Q(F-P) & ((F+P)^\top  Q)^\top   \\(F-P)^\top  Q  & -I_m & D^\top  \\ (F+P)^\top  Q & D & -I_m \rix\\
&=  \mat{ccc}Q &0& 0\\0 &I_m& 0\\ 0& 0&I_m\rix
\mat{ccc} -2R & (F-P) & (F+P)\\(F-P)^\top   & -I_m & D^\top  \\ (F+P)^\top   & D & -I_m  \rix
\mat{ccc}Q &0& 0\\0 &I_m& 0\\ 0& 0&I_m\rix \prec 0,
\end{align*}
since $Q$ is symmetric positive definite and $K_s \succ 0$.
\end{proof}

In order to prove that the converse of Theorem~\ref{thm:sphtoscat} is also true, we define the sPH form of a system~\eqref{eq1:sys}.

\begin{definition}
    A system $(A,B,C,D)$ is said to admit a sPH form if there exists a PH system as defined in~\eqref{eq:phsystem1} such that 
    $A=(J-R)Q$, $B=F-P$, $C=(F+P)^\top Q$ and $K_s \succeq 0$, where $K_s$ is defined in~\eqref{eq:bourealK}. 
\end{definition}

\begin{theorem} \label{thm:main_result} 
    Let $\Sigma=(A,B,C,D)$ be a system in the form~\eqref{eq1:sys}.
    The system $\Sigma$ admits an sPH form if there exists a symmetric positive definite matrix $X\in \R^{n,n}$ such that ${BR}_{LMI}(X) \preceq 0$, where ${BR}_{LMI}(X)$ be defined by~\eqref{eq:LMI1strict}.
\end{theorem}
\begin{proof}
Suppose that there exists $X \succ 0$ such that 
$\revise{{BR}_{LMI}(X)} \preceq 0$. Define
\begin{equation*}
 J:=\frac{AX^{-1}-(AX^{-1})^\top  }{2}, \quad R:=-\frac{AX^{-1}+(AX^{-1})^\top  }{2}, \quad Q:=X, 
 \vspace{-0.4cm}
\end{equation*}
\begin{equation*}
F:=\frac{1}{2}(B+X^{-1}C^\top  ), \quad \text{and}\quad P:=\frac{1}{2}(-B+X^{-1}C^\top  ).  \label{DHformcstr}
\end{equation*}
Clearly, $Q \succ 0$ and we have
\[
(J-R)Q=A,\quad F-P=B, \quad (F+P)^\top  Q =C.
\]
Further, we have 
\begin{eqnarray*}
K_s &=&\mat{ccc}2R & -(F-P) & -(F+P)\\ -(F-P)^\top   & I_m & -D^\top  \\ -(F+P)^\top   & -D & I_m \rix = -\mat{ccc} AX^{-1}+X^{-1}A^\top   & B & X^{-1}C^\top   \\ B^\top   & -I_m & D^\top  \\ CX^{-1} & D & -I_m \rix \\
&=& - \mat{ccc} X^{-1} & 0 & 0 \\0& I_m & 0 \\0 & 0& I_m \rix
\mat{ccc} A^\top  X +X A & XB & C^\top   \\B^\top  X &-I_m & D^\top  \\ C & D & -I_m \rix
\mat{ccc} X^{-1} & 0 & 0 \\0& I_m & 0 \\0 & 0& I_m \rix \succeq 0,
\end{eqnarray*}
since $X$ is a symmetric positive definite solution of the LMI ${BR}_{LMI}(X)\preceq 0$. 
This implies that $\Sigma$ admits a sPH form.
\end{proof}

\reviseELA{Let us summarize} the various equivalent conditions of a system to be scattering passive. 
\reviseELA{The first three equivalent statements were known previously, and we refer the interested reader to~\cite{BroLME20}. The fourth statement is our contribution, and follows from our Theorems~\ref{thm:sphtoscat} and~\ref{thm:main_result}.}

\begin{theorem} \label{thm:suff_BR}
Consider a minimal system $\Sigma=(A,B,C,D)$ in the form~\eqref{eq1:sys} and let \revise{${BR}_{LMI}(X)$} be defined by~\eqref{eq:LMI1strict}. Then the following are equivalent.
\begin{enumerate}
    \item  $\Sigma$ is BR.
    \item The \revise{LMI} ${BR}_{LMI}(X) \preceq 0$ has a symmetric positive definite solution $ X \in \R^{n,n}$.
    \item $\Sigma$ is scattering passive.
    \item $\Sigma$ admits a sPH form.
\end{enumerate}
\end{theorem}
\begin{proof}
    For the proof of \revise{the equivalence of} the first three statements, we refer to~\cite{BroLME20}. \\
    (2)$\implies$(4): this follows from Theorem~\ref{thm:main_result}.\\
    (4)$\implies$(2): the proof is identical to Theorem~\ref{thm:sphtoscat} as $X=Q$ becomes a solution to the matrix equation ${BR}_{LMI}(X) \preceq 0$.
\end{proof}

\begin{remark}{\rm
    Since BRness, scattering passivity, and the solution of the matrix equation 
    \mbox{${BR}_{LMI}(X)\preceq 0$} are all equivalent for minimal systems according to Theorem~\ref{thm:suff_BR}, they will be used interchangeably throughout the manuscript. Note that if we do not consider minimal representations, there is an intricate relationship between the different notions~\cite{CheGH23,CheGHM23}. Since we only consider numerically minimal realizations, this is not discussed further.
    }
\end{remark}

Theorem~\ref{thm:suff_BR} provides a new characterization of scattering passive systems in terms of sPH systems, which will be used in the next section to find an approximate solution to the nearest BR system problem~\eqref{probdef}.

\section{Formulation of the nearest BR system problem in an sPH form} 
\label{sec:NearsPH}

In this section, we exploit Theorem~\ref{thm:suff_BR} and propose a formulation to solve the nearest BR system problem~\eqref{probdef} by finding the nearby BR system in an sPH form. This formulation also leads to an algorithmic way to check if a system is BR.

Following Theorem~\ref{thm:suff_BR}, this problem can be reformulated as the following optimization problem 
\begin{align}  
\min_{Z \in \mathcal Z, J^\top  = -J, Q \succ 0} &  
\left\|
A-(J-Z_{11}/2) Q \right\|_F^2  +  
\left\|B + Z_{12}\right\|_F^2 + 
\big\|  C + Z_{13}^\top  Q \big\|_F^2  + \left\| D + Z_{23}^\top  \right\|_F^2,  \label{eq:mainoptprob}
\end{align}
where $Z$ is defined by blocks as follows 
\begin{equation} \label{eq:defZ}
Z = 
 \mat{ccc} 
Z_{11} & Z_{12} & Z_{13} \\ 
Z_{12}^\top   & I_m & Z_{23} \\ 
Z_{13}^\top   & Z_{23}^\top  & I_m 
\rix = 
\mat{ccc} 
2R & -(F-P) & -(F+P) \\ 
 -(F-P)^\top   & I_m & -\tilde{D}^\top  \\ 
 -(F+P)^\top   & -\tilde{D} & I_m
\rix \succeq 0. 
\end{equation}  
The matrix $Z$ encodes the positive semidefinite (PSD) constraint in the sPH form. The feasible set for $Z$ is 
\[ 
\mathcal Z = \left\{ Z \in \mathbb{R}^{n+2m,  n+2m} \ : \ Z = 
 \mat{ccc} 
Z_{11} & Z_{12} & Z_{13} \\ 
Z_{12}^\top   & I_m & Z_{23} \\ 
Z_{13}^\top   & Z_{23}^\top  & I_m  
\rix \succeq 0 \right\}. 
\] 

By Theorem~\ref{thm:suff_BR}, solving~\eqref{eq:mainoptprob} allows one to find the nearest BR system in an sPH form of a given system $(A,B,C,D)$. The distance here is measured as the Frobenius norm between the given system, $(A,B,C,D)$, and the sought BR system,  
$(\tilde A, \tilde B, \tilde C, \tilde D) = 
\left( (J-R)Q, F-P, (F+P)^\top Q, \tilde D\right)$, defined in the sPH 
form.  

   Other error measures could be used, e.g., $H_{2}$-norm or $H_{\infty}$-norm of the transfer function~\cite{burke2006hifoo, GuPBV12, SchV23, nicodemus2025klap}.  
\reviseELA{
In fact, our formulation~\eqref{eq:mainoptprob} depends on the representation $(A,B,C,D)$ of the system, not on the transfer function. Since this representation is not unique in general, the solution that our algorithm will compute by solving \eqref{eq:mainoptprob} will be sensitive to the representation. 
A possible way to control this sensitivity is to use} a weighted norm  
\begin{equation} \label{eq:weightedobjfun}
f^w(J,Z,Q) := 
w_1 \left\|A-(J-Z_{11}/2) Q \right\|_F^2  
+ w_2 \left\|B + Z_{12}\right\|_F^2 
+ w_3 \big\|  C + Z_{13}^\top  Q \big\|_F^2  
+ w_4 \left\| D + Z_{23}^\top  \right\|_F^2, 
\end{equation} 
for some nonnegative weights \revise{$w_i \geq 0$ for $i=1,\dots,4$}. We will design our algorithm for this weighted norm in Section~\ref{sec:EAOnearestBRsPHform}.  
Considering other error measures is a topic of further research.

Unfortunately,~\eqref{eq:mainoptprob} is a non-convex optimization problem, and hence we cannot hope in general to find a globally optimal solution. However, the feasible set in~\eqref{eq:mainoptprob} is a convex set that can be useful in optimization methods, see~Section~\ref{sec:eaocheckBR}. 

In Section~\ref{sec:checkBR}, we discuss the case when the input system is actually BR, in which case  \eqref{eq:mainoptprob} can be reformulated as an SDP. Even when the input system is not BR, the solution of this SDP can be used as an initialization for solving~\eqref{eq:mainoptprob}. For large-scale problems, solving this SDP might be impractical with standard software, and we propose a simple alternating optimization (AO) strategy for it. 

In Section~\ref{sec:EAOnearestBRsPHform}, we propose a similar AO strategy to solve \eqref{eq:mainoptprob}. We first rely on interior-point methods (IPMs) to solve the subproblems that can handle medium-scale problems ($n$ and $m$ up to a few tens on a standard laptop). To handle larger problems ($n$ and $m$ up to a few hundred on a standard laptop), we implement fast gradient methods to solve these subproblems.

%so that $Z_{11} \in \mathbb{R}^{n,  n}$, $Z_{22} \in \mathbb{R}^{m,  m}$, $Z_{33} \in \mathbb{R}^{m,  m}$. 
%Note that, in the impedance case, $N$ is fixed because there is no constraints beyond skew-symmetry \ngc{no, $N$ is involved in the PSD constraints within $Z$\dots}, and hence $N$ can be taken equal to the skew-symmetric part of $D$ that is, $N = (D-D^\top)/2$. Hence we can assume w.l.o.g.\ that $D$ is symmetric. In the impedance case, $N$ can be optimized as a variable; see~\cite{gillis2018finding}. 

\subsection{Checking BRness} \label{sec:checkBR}

If a system is BR, then there exists a system in sPH form such that the optimal objective function value of~\eqref{eq:mainoptprob} is equal to zero; see Theorem~\ref{thm:suff_BR}. 
In particular, there exists an optimal solution of~\eqref{eq:mainoptprob} such that $A = (J-R)Q$ and $C = (F+P)^\top  Q$ where $Q \succ 0$.   
Hence, similarly as done in~\cite{gillis2017computing}, we can make the following change of variable: $\hat Q := Q^{-1} \succ 0$ so that $A \hat Q = J-R$ and $C \hat Q = (F+P)^\top $, and solve instead 
\begin{equation}  
\min_{Z \in \mathcal Z, J^\top  = -J, \hat Q \succeq \epsilon I_n}   
\hat f\big(J,Z,\hat Q\big) := \left\|
A \hat Q -(J-Z_{11}/2) \right\|_F^2  +  
\left\|B + Z_{12}\right\|_F^2 + \left\|  C \hat Q + Z_{13}^\top  \right\|_F^2  + \left\| D \revise{+} Z_{23}^\top  \right\|_F^2.   \label{eq:mainoptprob2}
\end{equation} 
% where 
% \[
% f^i\big(Z,J,Q^i\big) = \left\|
% A Q^i -(J-Z_{11}/2) \right\|_F^2  +  
% \left\|B + Z_{12}\right\|_F^2 + \left\|  C Q^i + Z_{13}^\top \right\|_F^2  + \left\| D - Z_{32} \right\|_F^2. 
% \]
We use $\hat Q \succeq  \epsilon I_n$ for some $0 < \epsilon \ll 1$ sufficiently small, to avoid the trivial solution (that is, all variables equal to zero).   
Problem~\eqref{eq:mainoptprob2}
 is an SDP hence can be solved efficiently. 

%The system $(A,B,C,D)$ is BR if and only if there exist an optimal solution to~\eqref{eq:mainoptprob2} for $\epsilon$ sufficiently small  with objective function value equal to zero; this follows from Theorem~\ref{thm:suff_BR}. 
\revise{If the system $(A,B,C,D)$ is not BR,} the optimal solution of~\eqref{eq:mainoptprob2} provides an approximate solution to~\eqref{eq:mainoptprob} which could be used as an initialization for an algorithm that tackles the non-convex problem~\eqref{eq:mainoptprob}, similarly as done in~\cite{gillis2018finding} for PR systems. 

We make the following two observations: 
\begin{itemize}
\item The problem~\eqref{eq:mainoptprob2} does not necessarily have a unique solution, especially if $(A,B,C,D)$ is BR. This could possibly be understood by noting that the construction of the sPH form in Theorem~\ref{thm:main_result} is not unique, as is the solution $X \succ 0$ of the matrix equation $BR_{LMI}(X)\preceq 0$, where ${BR}_{LMI}(X)$ is as defined in~\eqref{eq:LMI1strict}.

    \item To reduce the number of variables in~\eqref{eq:mainoptprob2}, observe that the optimal $J$ is given by the skew symmetric part of $A\hat Q$, which we denote $\skewp(A\hat Q)$. 
    Hence, one can get rid of variable $J$ and replace the first term by 
    $\left\|
\sym(A \hat Q) + Z_{11}/2 \right\|_F^2$,  where $\sym(A) = (A+A^\top )/2$ is the symmetric part of matrix $A$.  
\end{itemize} 

To solve the convex optimization problem~\eqref{eq:mainoptprob2}, we propose two strategies described in the next two subsections. 
\revise{
The first one uses interior-point methods, hence has a quadratic rate of convergence but is computationally expensive. The second-one uses a first-order method with typical linear convergence, hence is computationally much cheaper but may converge slowly. 
}

\subsubsection{Interior-point methods} \label{sec:checkBRIPM}

One can use a SDP solver to solve \eqref{eq:mainoptprob2} directly, 
and we use CVX modeling system~\cite{cvx, gb08} with the IPM solver SDPT3 4.0~\cite{toh1999sdpt3, tutuncu2003solving}. 
See the code \texttt{nearest$\_$inv$\_$BRsPHform.m} with \texttt{options.algo = `CVX-IPM'}.  
This code can be used for matrices of the order $n+2m$ up to about a hundred on a standard laptop; see Section~\ref{sec:numexp} for numerical experiments. 
For example, we will use it on a system with $n=200$ and $m=1$ in Section~\ref{sec:karim}, and it takes about 25 minutes to obtain an optimal solution.

\subsubsection{Extrapolated Alternating Optimization (E-AO)}  \label{sec:eaocheckBR}

To address larger instances, we propose to solve~\eqref{eq:mainoptprob2} using AO: we alternatively update the variables $\hat Q$ and $(J,Z)$. 

\paragraph{Optimizing $\hat Q$} The subproblem in $\hat Q$, for $(J,Z)$ fixed, is a PSD Procrustes problem. Using a fast gradient method (FGM)~\cite{nesterov2013introductory}, an optimal first-order method, it requires $O(n^3)$ operations per iteration to compute the gradient and project onto the feasible set. In fact, the projection onto the set of PSD matrices requires an eigendecomposition and then simply clipping the negative eigenvalues to zero~\cite{higham1988computing}. 
We use the code developed in~\cite{gillis2018semi} (\texttt{Procrustes$\_$FGM.m}).  
In practice, we will run FGM for only a few iterations, as it is not necessary to solve the subproblem with high accuracy since $(J,Z)$ will be modified at the next iteration. 

\paragraph{Optimizing $(J,Z)$} As explained above, the optimal $J$ is given by $\skewp(A\hat Q)$. %while the problem in $Z$ is a quadratic problem with the Hessian being the identity matrix. 
For $Z$, given $\hat Q$, we need to solve  
\[ 
\min_{Z \in \mathcal Z} 
\left\|
\sym(A \hat Q) + Z_{11}/2 \right\|_F^2  +  
\left\|B + Z_{12}\right\|_F^2 + \left\|  \hat Q C^\top   + Z_{13} \right\|_F^2  + \left\| D + Z_{23}^\top  \right\|_F^2, 
\] 
which is equivalent to 
\begin{equation}
    \min_{Z \in \mathcal Z} 
 \left\|
Z_{11} - (-2 \sym(A \hat Q))  \right\|_F^2  +  
4 \left\| Z_{12} - (-B) \right\|_F^2 + 
4 \left\|  Z_{13} - (-\hat Q C^\top )  \right\|_F^2  + 4 \left\|  Z_{23}^\top  - (-D)\right\|_F^2. \label{eq:problemZcheck}
\end{equation}
This is also a PSD constrained least squares problem, like for $\hat Q$, that can be solved with FGM. The only tricky part is the projection onto the feasible set, $\mathcal Z$, which is described below. 

To initialize $\mathcal Z$, we will solve 
\begin{equation}
\min_{Z \in \mathcal Z} \left\| 
 \mat{ccc} 
Z_{11} & Z_{12} & Z_{13} \\ 
Z_{12}^\top   & I_m & Z_{23} \\ 
Z_{13}^\top   & Z_{23}^\top  & I_m 
\rix
- 
 \mat{ccc} 
- 2 \sym(A \hat Q) & -B & - \hat Q C^\top   \\ 
-B^\top   & I_m & -D^\top  \\ 
- C \hat Q  & -D & I_m 
\rix
\right\|_F^2, \label{eq:projectionZ}
\end{equation}
which is the projection of a matrix onto the set $\mathcal Z$. 
Note that \eqref{eq:projectionZ} is, unfortunately, not equivalent to~\eqref{eq:problemZcheck}; in fact, \eqref{eq:projectionZ} is equivalent to 
\[ 
    \min_{Z \in \mathcal Z} 
  \left\|
Z_{11} - (-2 \sym(A \hat Q))  \right\|_F^2  +  
2 \left\| Z_{12} - (-B) \right\|_F^2 + 
2 \left\|  Z_{13} - (-\hat Q C^\top )  \right\|_F^2  
+ 2 \left\|  Z_{23}^\top  - (-D)\right\|_F^2. 
\] 

\paragraph{Projecting onto $\mathcal Z$} Projecting onto $\mathcal Z$ requires to solve another SDP. However, we resort to a standard  strategy to project onto $\mathcal Z$, namely the alternating direction method of multipliers (ADMM)~\cite{boyd2011distributed}, 
which is reminiscent of alternating projection (AP) introduced by von Neumann in 1933 (see, e.g.,~\cite{bauschke1996projection} and the references therein). 
\revise{For simplicity, let us denote $Z \notin \mathcal Z$ an iterate obtained after a gradient step of FGM applied on~\eqref{eq:problemZcheck}}. 
We need to project it back onto $\mathcal Z$, that is, solve $\min_{X \in \mathcal Z} \|X - Z\|_F^2$ where $\mathcal Z =   \mathcal S_+ \cap \mathcal{D}$,  $\mathcal S_+$ is the (convex) set of PSD matrices of appropriate dimension, and $\mathcal{D}$ is the affine set 
\[
\mathcal{D} = \left\{ Z \in \mathbb{R}^{n+2m, n+2m} \ : \ 
Z = \mat{ccc} 
Z_{11} & Z_{12} & Z_{13} \\ 
Z_{12}^\top   & I_m & Z_{23} \\ 
Z_{13}^\top   & Z_{23}^\top  & I_m  
\rix \right\}. 
\]  
ADMM first duplicates the variable $X=Y$, and considers the equivalent problem 
\[
\min_{X \in \mathcal S_+, Y \in \mathcal{D}} 
 \|X-Z\|_F^2 
+  \|Y-Z\|_F^2 \quad \text{ such that } \quad X = Y. 
\]
Then it introduces the Lagrange multpliers $\Lambda \in \mathbb{R}^{n+2m,n+2m}$, and considers the augmented Lagrangian 
\[
\min_{X \in \mathcal S_+, Y \in \mathcal{D}} 
 \|X-Z\|_F^2 
+  \|Y-Z\|_F^2 
- 2 \rho \left\langle \Lambda, X-Y \right\rangle
+ \rho  \|X-Y\|_F^2, 
\] 
which can be rewritten as 
\[
\min_{X \in \mathcal S_+, Y \in \mathcal{D}} 
 \|X-Z\|_F^2 
+  \|Y-Z\|_F^2 
+ \rho \|X-Y-\Lambda\|_F^2. 
\] 
Finally ADMM alternates between the update of $X$, $Y$ 
and $\Lambda$. The updates of $X$ and $Y$ reduce to projections on their respective feasible set, while the update of $\Lambda$ penalizes the difference between $X$ and $Y$: letting $\lambda = \frac{1}{1+\rho}$, 
and initializing $\Lambda=0$, $Y=Z$, ADMM updates $X$, $Y$ and $\Lambda$ as follows. 
\begin{equation*}
 \text{ ADMM to project onto $\mathcal{Z}$:} 
 \quad \left\{ 
 \begin{array}{lll}
      X & \leftarrow & \mathcal P_{\mathcal S_+} \big( \lambda Z + (1-\lambda)  (Y+\revise{\Lambda})  \big), \\ 
        Y & \leftarrow &  \mathcal P_{\mathcal D} \big( \lambda Z + (1-\lambda) (X-\revise{\Lambda})  \big), \\
            \Lambda & \leftarrow &  \Lambda-(X-Y), 
            \end{array} 
            \right. 
\end{equation*}
where $\mathcal P_{\mathcal S}(\cdot)$ is the projection onto the set $\mathcal S$. 
The projection onto $\mathcal{D}$ is trivial, simply replacing the (2,2) and (3,3) blocks by the identity matrix. 
Projection onto $\mathcal S_+$ 
costs $O\big((n+2m)^3\big)$ operations as it requires an eigendecomposition.  
We refer the interested reader to~\cite{bilkova2021projection} for more details on ADMM used for projecting on the intersection of two sets.   

Let 
\[
X = \mat{ccc} 
X_{11} & X_{12} & X_{13} \\ 
X_{12}^\top   & X_{22} & X_{23} \\ 
X_{13}^\top   & X_{23}^\top  & X_{33}
\rix \succeq 0, 
\]
be the last iterate generated by ADMM. We typically have $X_{22} \neq I_m$ or $X_{33} \neq I_m$, although $X_{22} \approx I_m$ or $X_{33} \approx I_m$ if ADMM is run for sufficiently many iterations since it converges to the projection of $Z$ onto $\mathcal{Z}$. 
In order to make $X$ feasible, let $\alpha_i = \max\left(1,\lambda_{\max}(X_{ii})\right)$ for $i=2,3$, where $\lambda_{\max}(A)$ denotes the largest eigenvalue of  $A$, 
which implies that $\alpha_i I \succeq X_{ii}$ and that  
\[
\mat{ccc} 
X_{11} & X_{12} & X_{13} \\ 
X_{12}^\top   & \alpha_2 I_m & X_{23} \\ 
X_{13}^\top   & X_{23}^\top  & \alpha_3 I_m 
\rix  \succeq X \succeq 0. 
\] 
Hence, using a scaling by $\alpha_2$ and $\alpha_3$, 
\begin{equation} \label{eq:Zclip}
\tilde X = 
\mat{ccc} 
X_{11}                & X_{12}/\alpha_2                   & X_{13}/\alpha_3 \\ 
X_{12}^\top /\alpha_2  &  I_m                                & X_{23}/(\alpha_2 \alpha_3) \\ 
X_{13}^\top /\alpha_3  & X_{23}^\top /(\alpha_2 \alpha_3)  & I_m  
 \rix 
 %\succeq
% \mat{ccc} 
% X_{11}/\alpha & X_{12}/\alpha & X_{13}/\alpha \\ 
% X_{12}^\top/\alpha  &  I & X_{23}/\alpha \\ 
% X_{13}^\top/\alpha  & X_{23}^\top/\alpha & I 
% \rix 
\in \mathcal{Z}.  
\end{equation}
If $X_{22}$ and $X_{33}$ are close to the identity, $\alpha_2, \alpha_3 \approx 1$ and $\tilde X \approx X$. 
This can be controlled in practice by monitoring the convergence of ADMM via the quantity 
$\left \| X-Y \right\|_F$. 

% Can we prove some error bounds on this scaling, assuming $Z$ is close to the set $\mathcal{D}$, that is, assuming $\|Z_{ii}-I\|_F \leq \epsilon$ for $i=1,2$ and some small $\epsilon$ (maybe in relative error? $\|Z_{ii}-I\|_F \leq \epsilon \|Z_{ii}\|_F$?) 
% Assume 
% \[
% \min \left( \lambda_{\min}(Z_{22}), \lambda_{\min}(Z_{3})\right) \geq 1-\epsilon , 
% \] 
% then $\alpha_2, \alpha_3 \leq \frac{1}{1-\epsilon}$. 
% In particular, 
% if $\|Z-I\|_F \leq \epsilon$, 
% then $\lambda_{\min}(Z)$ 
% Can we do better than this simple scaling? 
%We will use $ \rho = 1$. Note that RAP with $\rho \sigma_k = 1$ for all $k$ reduces to AP. 
% Figure~\ref{fig:accelRAP} illustrates the acceleration effect of RAP on randomly generated instances: We generate $Z = WW^\top \succeq 0$ where the entries of $W$ are sampled with the normal distribution, that is, 
% \texttt{W = randn(n+2*m,n+2*m)}, and\footnote{It turns out generating $Z$ in this way make it have $30$ negative eigenvalues with very high probability: it was the case for 1000 matrices we generated in this way.} then let $Z_{22} = Z_{33} = I_m$. We used $n=25$, $m=15$. 
% \begin{figure}[ht!]
% \begin{center}
% \includegraphics[width=7cm]{figures/accelRAP.png} 
% \caption{Convergence of 
% AP and RAP on 10 randomly generated 55-by-55 matrices.}
% \label{fig:accelRAP}
% \end{center}
% \end{figure}
%We will use $\rho = 0.25$ in the rest of the paper. 

\paragraph{Acceleration with extrapolation} 

To accelerate the convergence of AO, we use extrapolation for block coordinate descent methods; see~\cite{HPG_TITAN} and the references therein. 
In a nutshell, given two consecutive iterates, $x^k$ and $x^{k-1}$, an extrapolated point is constructed as $\hat x^{k} = x^k + \beta_k ( x^k - x^{k-1} )$, 
and the next iterate, $x^{k+1}$, is computed from this extrapolated point. 
Intuitively, $\hat x^{k}$ pushes $x^k$ in the direction defined by the two previous iterates. In convex optimization, extrapolation leads to acceleration~\cite{nesterov2017efficiency}. In non-convex optimization, this scheme can guarantee convergence under some suitable choices of $\beta_k$, but does not guarantee acceleration in general (it does in some specific cases, e.g., for matrix factorization~\cite{xuprovable24}). 
For simplicity, we resort to a heuristic extrapolation with restart, similar to what is proposed in~\cite{o2015adaptive}; 
see Algorithm~\ref{alg:checkBRfgm}. The restarting procedure guarantees that the objective value goes down at each iteration. 
We illustrate the empirical acceleration effect of extrapolation in Section~\ref{sec:numexp}.  
\reviseELA{Because of the heuristic extrapolation choice, and the approximate projection on the feasible set, we cannot prove the convergence of the iterates for Algorithm~\ref{alg:checkBRfgm}. This is a rather technical topic, see, e.g., \cite{golbabaee2018inexact}, and left for further research.}

\begin{algorithm}[ht!] 
\caption{Extrapolated Alternating Optimization (E-AO) to check BRness, that is, solve~\eqref{eq:mainoptprob2}} \label{alg:checkBRfgm}
\begin{algorithmic}[1]

\REQUIRE System $(A,B,C,D)$, 
initial point $\hat Q_0 \succ 0$, 
maximum number of iteration\revise{s} $K$, 
number of inner iterations $I$, 
ADMM parameter $\rho$,  
stopping criterion $\delta \ll 1$, 
extrapolation parameters $\beta \in [0,1)$,  
$\eta > \gamma > 1$, \revise{lower bound $0 < \epsilon \ll 1$ on the eigenvalues of $\hat Q$}. 

(Default: $\hat Q_0 = I_n$, 
$K=1000$, 
$I=2$, 
$\delta = 10^{-6}$, 
$\beta = 0.5, 
\eta = 2, 
\gamma = 1.05$, 
$\rho = 10$, \revise{$\epsilon = 10^{-4}$}.)

\ENSURE An approximate solution $(Z, J, \hat Q)$ to~\eqref{eq:mainoptprob2}. 

\medskip 

\STATE $\hat Q_e \revise{ = \hat Q} = \hat Q_0$, $k = 1$. 

\WHILE{
$k \leq K$ 
and $\left( k \leq 3 \text{ or } e_{k-2} - e_{k-1}  \geq \delta e_{1} \right)$ 
}
%\Statex \Comment {\textbf{Step $1$: Update $u$}}
 \STATE  \textbf{\emph{\% Update $(Z,J)$}} 

   \STATE \textbf{if} $k \geq 2$ 
  \textbf{then} $(Z_p,J_p) = (Z,J)$. 
  \textbf{end if} 
  \emph{\% Keep the previous iterate in memory}  
  
 \STATE Compute $Z$ using $I$ iterations of FGM on 
 \[
    \min_{Z \in \mathcal Z} 
 \left\|
Z_{11} - (-2 \sym(A \hat Q_e))  \right\|_F^2  +  
4 \left\| Z_{12} +B \right\|_F^2 + 
4 \left\|  Z_{13} + \hat Q_e C^\top   \right\|_F^2  + 4 \left\|  Z_{23}^\top  + D\right\|_F^2. 
 \] 
 The projection is performed with ADMM with parameter $\rho$. 

 \STATE Let $J = \skewp(A \hat Q_e)$. 

 \IF {$k=1$}   

\STATE $(Z_e,J_e) = (Z,J)$. 

\ELSE 

\STATE $(Z_e,J_e) = (Z,J) + \beta \big((Z,J)-(Z_p,J_p)\big)$. \emph{\% Extrapolation} 

\ENDIF

 %\Statex \Comment {\textbf{Step $2$: Update $v$}}
 
 \STATE \textbf{\emph{\% Update \revise{$\hat Q$}}}
  \STATE $\hat Q_p = \hat Q$. \emph{\% Keep the previous iterate in memory}  
  
 \STATE Update $\hat Q$ using $I$ iterations of FGM on 
 \[
\min_{\hat Q \succeq \epsilon I_n} \| A  \hat Q  - (J_e-(Z_e)_{11}/2) ||_F^2 + \| C \hat Q + (Z_e)_{13}^\top  \|_F^2 . 
 \]

\STATE $\hat Q_e = \hat Q + \beta (\hat Q-\hat Q_p).$ \emph{\% Extrapolation}

\STATE $e_k \leftarrow \sqrt{\frac{\hat f(J,Z,\hat Q)}{\hat f(0,0,0)}}$. \emph{\% Relative error}  

 \STATE  \textbf{\emph{\% Restart scheme when the objective increases}} 
 
\IF {\revise{$k \geq 2$ and} $e_k > e_{k-1}$ }   

\STATE $(Z,J) = (Z_p,J_p)$, \revise{$\hat Q_p = \hat Q$}, \reviseELA{$e_k = e_{k-1}$}, $\beta = \frac{\beta}{\eta}, 
I = \lceil 1.1 \ I \rceil$. \emph{\% Restart}  

\ELSE 

\STATE $\beta \leftarrow \min( 1, \gamma  \beta).$ 

\ENDIF

\STATE $k \leftarrow k+1$. 

\ENDWHILE  

\end{algorithmic}
\end{algorithm}

\subsection{E-AO algorithm for the nearest BR problem~\eqref{eq:mainoptprob}}  \label{sec:EAOnearestBRsPHform}

In this section, we propose an algorithm to tackle a weighted version of~\eqref{eq:mainoptprob}, namely 
\begin{equation} \label{eq:mainoptprobweight}
\min_{Z \in \mathcal Z, J^\top = -J,Q \succeq 0} f^w(J,Z,Q), 
\end{equation} 
where $f^w(J,Z,Q)$ is as defined in~\eqref{eq:weightedobjfun}, 
in order to find the nearest BR system in an sPH form. The proposed algorithm follows the same strategy as Algorithm~\ref{alg:checkBRfgm}, namely an E-AO approach. 
Although~\eqref{eq:mainoptprob} is not globally convex, it is bi-convex in variables $Q$ vs.\ $(J,Z)$. In fact, the subproblems in $Q$ for $(J,Z)$ fixed, and vice versa, are SDPs.  
It is therefore natural to use an AO strategy by optimizing $Q$ and $(J,Z)$ alternatively. 
%This is the same strategy as for solving~\eqref{eq:mainoptprob2} in the previous section, but the subproblem in $(J,Z)$ do not have closed form anymore. 
Algorithm~\ref{alg:eao} summarizes the proposed algorithm, which has a very similar flavor to Algorithm~\ref{alg:checkBRfgm} and also incorporates extrapolation whose acceleration effect will be illustrated in Section~\ref{sec:numexp}. 

%The main difference is that $(J,Z)$ cannot be updated in closed form but their update requires to solve an SDP.  

\begin{algorithm}[h!]
\caption{Extrapolated Alternating Optimization (E-AO) to find the nearest BR system in an sPH form, that is, solve~\eqref{eq:mainoptprob}} \label{alg:eao}
\begin{algorithmic}[1]

\REQUIRE System $(A,B,C,D)$, 
initial point $(J,Z,Q)$, 
maximum number of iterations $K$, 
stopping criterion $\delta \ll 1$, 
ADMM parameter $\rho$, 
extrapolation parameters $\beta \in [0,1)$,  
$\eta > \gamma > 1$, weight vector $w \in \mathbb{R}^{4}_+$ defining the weighted objective $f^w(J,Z,Q)$~\eqref{eq:weightedobjfun} 

(Default: algo = FGM, 
$Q = I_n$, 
$J = \skewp(A)$, 
$Z$ in \eqref{eq:projectionZ} solved with ADMM, 
$K=1000$, 
$I=10$, 
$\delta = 10^{-6}$, 
$\beta = 0.5, 
\eta = 2, 
\gamma = 1.05$, 
$\rho = 10$, $w = [1,1,1,1]$.)

\ENSURE An approximate solution to~\eqref{eq:mainoptprob}. 

\medskip 

\STATE $Q_e = Q$, $k = 1$. 

\WHILE{
$k \leq K$ 
and $\left( k \leq 4 \text{ or } e_{k-2} - e_{k-1}  \geq \delta e_{1} \right)$ 
}
%\Statex \Comment {\textbf{Step $1$: Update $u$}}
 \STATE  \textbf{\emph{\% Update $(Z,J)$}} 

  \STATE  $(Z_p,J_p) = (Z,J)$. \emph{\% Keep the previous iterate in memory}  

  \IF{algo = `IPM'}
  
 \STATE  Solve~\eqref{eq:mainoptprob} in variables $(Z,J)$  for $Q=Q_e$ fixed using CVX. 
 
 \ELSIF{algo = `FGM'}
 
  \STATE From $(J,Z)$, apply $I$ iterations of FGM on 
  the problem $\min_{Z \in \mathcal Z, J^\top  = -J} f^w(J,Z,Q_e)$ 
  to obtain the new $(Z,J)$. 
  The projection of $Z$ onto $\mathcal{Z}$ is done approximately with $I$ iterations of ADMM \revise{with} parameter $\rho$.  
  
  \ENDIF 

%  \IF {$k=1$}   

% \STATE $(Z_e,J_e) = (Z,J)$. 

% \ELSE 

\STATE $(Z_e,J_e) = (Z,J) + \beta \big((Z,J)-(Z_p,J_p)\big)$. \emph{\% Extrapolation} 

%\ENDIF

 %\Statex \Comment {\textbf{Step $2$: Update $v$}}
 
 \STATE \textbf{\emph{\% Update $Q$}}
  \STATE $Q_p = Q$. \emph{\% Keep the previous iterate in memory}

  \IF{algo = `IPM'}
  
 \STATE  Solve~\eqref{eq:mainoptprob} in variables $Q$  for fixed $(J,Z) = (J_e,Z_e)$  using CVX and IPM.  
 
 \ELSIF{algo = `FGM'}
 
  \STATE Apply $I$ iterations of FGM on~\eqref{eq:mainoptprob} in variables $Q$ for $(Z,J) = (J_e,Z_e)$ fixed. 
  
  \ENDIF 

\STATE $Q_e = Q + \beta (Q-Q_p).$ \emph{\% Extrapolation}

\STATE $e_k \leftarrow \sqrt{\frac{f^w(J,Z,Q)}{f^w(0,0,0)}}$. \emph{\% Relative error}  

 \STATE  \textbf{\emph{\% Restart scheme when the objective increases}} 
 
\IF {$k \geq 2$ and $e_k > e_{k-1}$}

\STATE $(Z,J) = (Z_p,J_p)$, $Q = Q_p$, \reviseELA{$e_k = e_{k-1}$}, $\beta = \frac{\beta}{\eta}$, $I = \lceil 1.1 \ I \rceil$. 

\ELSE 

\STATE $\beta \leftarrow \min( 1, \gamma  \beta).$ 

\ENDIF

\STATE $k \leftarrow k+1$. 

\ENDWHILE  

\end{algorithmic}
\end{algorithm}

The subproblems in variables $(J,Z)$ and $Q$ are SDPs. Again, we propose to solve them with two strategies: 
\begin{itemize}
    \item CVX and IPMs. The downside, again, is that the subproblems to solve are SDPs. Most general-purpose SDP solvers rely on second-order methods, scaling in $O((n+m)^6)$ operations, and hence we cannot handle systems with $n$ and $m$ larger than a few tens. 

\item FGM. With this method, we are able to handle larger systems. Note that we do not have an exact projection onto the set $\mathcal Z$ and resort to the use of ADMM, as in the previous section.  
    
\end{itemize}

We close this section with the following remark on systems that are both positive-real and bounded-real. 
\begin{remark}{\rm
In view of the definitions of impedance PH systems and scattering PH systems, we have that a PH system in the form~\eqref{eq:phsystem1} with $K_i \succ 0$ and $K_s \succ 0$, where $K_i$ and $K_s$ are respectively defined by~\eqref{eq:posrealK} and~\eqref{eq:bourealK}, is both positive-real (\cite[Theorem 3.2] {gillis2018finding}) and bounded-real (Theorem~\ref{thm:sphtoscat}). Thus, for a given system $\Sigma=(A,B,C,D)$, a nearby system which is both PR and BR can be found by finding a nearby PH system with $K_i \succ 0$ and $K_s \succ 0$. More precisely, solve the following optimization problem.
\begin{align}  
\min_{Z \in \mathcal Z, J^\top  = -J, Q \succ 0} &  
\left\|
A-(J-Z_{11}/2) Q \right\|_F^2  +  
\left\|B + Z_{12}\right\|_F^2 + 
\big\|  C + Z_{13}^\top  Q \big\|_F^2  + \left\| D + Z_{23}^\top  \right\|_F^2,  \label{eq:mainoptprobBRPR}
\end{align}
where $Z$ is a block matrix that encodes the positive semidefinite (PSD) constraints in both the iPH form and the sPH form. The feasible set for $Z$ is 
\[ 
\mathcal Z = \left\{ Z \in \mathbb{R}^{n+2m,  n+2m} \ : \ Z = 
 \mat{ccc} 
Z_{11} & Z_{12} & Z_{13} \\ 
Z_{12}^\top   & I_m & Z_{23} \\ 
Z_{13}^\top   & Z_{23}^\top  & I_m  
\rix \succeq 0,~\mat{cc} Z_{11} & Z_{12}-Z_{13}\\Z_{12}^\top - Z_{13}^\top & Z_{23}+Z_{23}^\top\rix \succeq 0
\right\}. 
\] 
We note that a method similar to Algorithm~\ref{alg:eao} can be used to solve~\eqref{eq:mainoptprobBRPR}.
}
\end{remark}

\revise{
\begin{remark}[Checking BRness with Algorithm~\ref{alg:eao}] \rm{
   Algorithm~\ref{alg:eao} can be used to check BRness, as done by Algorithm~\ref{alg:checkBRfgm}, since a solution with zero error would imply that the input system is BR. However,  Algorithm~\ref{alg:eao} tackles a non-convex optimization problem, and hence can be stuck on local optimal solutions. On the contrary, Algorithm~\ref{alg:checkBRfgm} solves a convex optimization problem~\eqref{eq:mainoptprob2}. 
   Of course, because Algorithm~\ref{alg:checkBRfgm} is a first-order method, it may converge slowly (see numerical examples in the next section) and hence it might be difficult to obtain solutions with high accuracy. 
   Also, since we are using an aggressive extrapolation strategy and an approximate projection, the proof of convergence is nontrivial, and is left for future work. We have however observed in practice a significantly better performance than a standard gradient descent; see the numerical experiments in Section~\ref{sec:numexp}. 
   If checking BRness is crucial, one should prefer to use the interior-point method (see Section~\ref{sec:checkBRIPM}) that provides guaranteed high-accuracy solutions. 
   }
\end{remark}
}

\reviseELA{
\begin{remark}[Convergence of Algorithm~\ref{alg:eao}]
    Because Algorithm~\ref{alg:eao} uses a restarting scheme that guarantees the objective to decrease, and since the objective is bounded below, the objective function values will converge. 
    However, convergence of the iterates is much more tricky, as we are in a non-convex setting, and we are using extrapolation rather aggressively (which works well in practice but does not come with convergence guarantees). Studying the convergence of iterates, possibly by adapting the algorithm, is a topic of further research.  
\end{remark}}
% \paragraph{Initialization} 

% \begin{itemize}
%     \item Identity. Since $Q=I$, the optimal solution of $(J,Z)$ can be computed... 

%     \item With the solution of~\eqref{eq:mainoptprob2} using CVX and IPM. 

%     \item With the solution of Algorithm~\ref{alg:checkBRfgm}. 
% \end{itemize}

\section{Numerical Experiments} \label{sec:numexp} 

In this section, we consider two problems. First, the SDP \eqref{eq:mainoptprob2}, which can be used to check whether a given system is BR. We will solve this problem with CVX and with Algorithm~\ref{alg:checkBRfgm}.  Then, we consider our main problem, finding the nearest BR system in an sPH form, that is, solving~\eqref{eq:mainoptprob}. To do so, we use E-AO, Algorithm~\ref{alg:eao}, which has two variants: EAO-IPM which solves the subproblems with CVX and IPMs, and EAO-FGM which solves the subproblems with a fast gradient method. 
To initialize Algorithm~\ref{alg:eao}, we will consider 3 approaches: 
\begin{itemize}
    \item Id: $Q = I_n$, in which case $J$ can be computed in closed form since $\hat Q = Q^{-1} = I_n$, and $Z$ is set as the projection in \eqref{eq:projectionZ}, as described in Section~\ref{sec:eaocheckBR} . 

    \item CVX: $(Q,J,Z)$ as the solution of~\eqref{eq:mainoptprob2} using CVX.  

    \item Alg1: $(Q,J,Z)$  as the solution of~\eqref{eq:mainoptprob2} solved with Algorithm~\ref{alg:checkBRfgm}. 
    
\end{itemize} 
We will denote EAO-alg(init) where alg $\in \{$IPM, FGM$\}$ and init $\in \{$Id, CVX, Alg1$\}$. For example EAO-IPM(Id) means  Algorithm~\ref{alg:eao} where subproblems are solved with CVX and it is initialized with $Q=I_n$ and $(J,Z)$ as described above. 

We will use the relative error to assess the quality of solutions, 
which is defined as $\sqrt{\frac{f^w(J,Z,Q)}{f^w(0,0,0)}}$ for the nearest BR system problem (see Algorithm~\ref{alg:eao}), and as $\sqrt{\frac{\hat f(J,Z,\hat Q)}{\hat f(0,0,0)}}$ for the problem of checking BRness 
 (see Algorithm~\ref{alg:checkBRfgm}). 
 Relative error is nicely interpretable as it reports the error relative to the norm of the system.  
 
 %in particular  $\sqrt{f^w(0,0,0)} = \sqrt{\hat f(0,0,0)} = \|(A,B,C,D)\|_F$ when the weights  in the objective are all ones. 

The numerical experiments were performed on a laptop with 
processor 12th Gen Intel(R) Core(TM) i9-12900H and 32,0Go RAM. 
The code, with which all numerical experiments can be performed directly, is available from~\url{https://gitlab.com/ngillis/nearestBRsysPHform}.

\subsection{Small system from~\cite{BoyBK89}} \label{kresex}

Let us start with a small system to illustrate the proposed models and algorithms. 
Consider the following standard LTI system  from~\cite[Section 6]{BoyBK89} 
\begin{equation}
A \text{$=$}\mat{cccc}-0.08 & 0.83& 0& 0\\-0.83& -0.08& 0& 0\\0 &0 &-0.7& 9\\0& 0& -9 &-0.7\rix, 
B \text{$=$} \mat{cc}
1 & 1 \\ 
0 & 0 \\
1 & -1 \\ 
0 & 0 \rix, 
C \text{$=$} \mat{cccc}0.4&0&0.4&0\\0.6&0&1&0\rix, 
D \text{$=$} \mat{cc}0.3 &0\\0 &-0.15 
\rix. \label{sys:boyd} 
\end{equation}

\paragraph{Checking BRness} 

This system is asymptotically stable but not BR. In fact, solving~\eqref{eq:mainoptprob2} with CVX provides a relative error of 4.31\% in 0.7 seconds, which means that there does not exist a solution $X \succ 0$ of the LMI $BR_{LMI}(X)\preceq 0$; see Theorem~\ref{thm:suff_BR}. Algorithm~\ref{alg:checkBRfgm} finds a solution with error 4.32\% within 0.12 seconds.  

%Figure~\ref{} shows the evolution of the error of Algorithm~\ref{alg:checkBRfgm}, with extrapolation ($\beta = 0.5$) and without extrapolation ($\beta = 0$). 

\paragraph{Nearest BR system in an sPH form}  

\begin{figure}[ht!]
\begin{center}
 \includegraphics[width=0.6\textwidth]{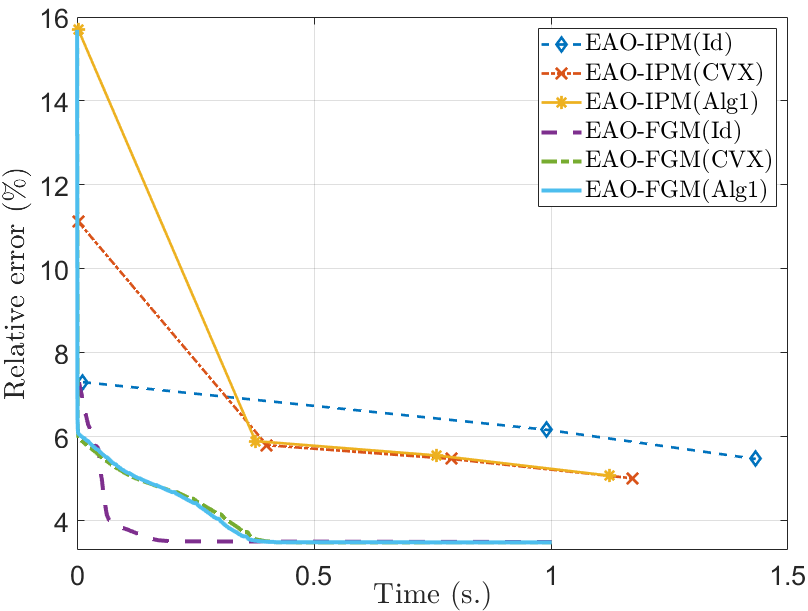}  
\caption{Relative error of the different variants of Algorithm~\ref{alg:eao} with different initializations for the system~\eqref{sys:boyd}. 
%Note that the curve $(\beta=0.5,I=2)$ on the left corresponds to the curve $\rho = 10$ on the right. 
\label{fig:ex1}}  
\end{center}
\end{figure}

Using Algorithm~\ref{alg:eao} to solve~\eqref{eq:mainoptprob} for this system, 
we obtain the following best approximation  with EAO-FGM (up to 3 digits of accuracy) 
\[
\tilde A = 
\left( \begin{array}{cccc} 
 -0.291 &  0.832 &  0.000 &  -0.002 \\ 
 -0.835 &  -0.267 &  0.000 &  -0.000 \\ 
 -0.003 &  0.003 &  -0.751 &  9.000 \\ 
 -0.030 &  0.005 &  -8.999 &  -0.751 \\ 
\end{array} \right) , 
\tilde B=\left( \begin{array}{cc} 
 0.928 &  0.945 \\ 
 -0.019 &  -0.019 \\ 
 0.961 &  -0.956 \\ 
 -0.015 &  -0.007 \\ 
\end{array} \right) ,~
    \]
 \[
\tilde C=\left( \begin{array}{cccc} 
 0.215 &  0.026 &  0.364 &  0.003 \\ 
 0.390 &  0.045 &  0.934 &  0.007 \\ 
\end{array} \right) ,~
\tilde D = \left( \begin{array}{cc} 
 0.213 &  -0.015 \\ 
 -0.115 &  -0.143 \\ 
\end{array} \right). 
\] 
Figure~\ref{fig:ex1} displays the evolution of the relative error, but does not account for the time to compute the initializations.  
EAO-FGM converges in less than half a second for all initializations, performing about 350 iterations, while EAO-IPM can only perform 3 to 4 iterations within 1 second, and does have time to converge.  EAO-FGM converges to the same solution for the 3 initialization strategies. 
This gives a nearby BR system with relative error 3.48\%. 
In terms of relative error for each matrix, we have 
\begin{equation*}
\frac{{\|A-\tilde A\|}_F}{{\|A\|}_F} = 2.29\%,~
\frac{{\|B-\tilde B\|}_F}{{\|B\|}_F} = 5.62\%,~
\frac{{\|C-\tilde C\|}_F}{{\|C\|}_F} = 22.69\%,~
\frac{{\|D-\tilde D\|}_F}{{\|D\|}_F} = 43.21\%.
\end{equation*} 
Note that the approximation error for $D$, and to a lesser extent for $C$,  is rather large because the norm of $D$ compared to the other matrices (in particular $A$) is smaller hence it implicitly has less importance in the objective function. 
By playing with the weight vector $w$ in the weighted objective function~\eqref{eq:weightedobjfun}, 
namely using $w=[0.5,2,5,20]$, we get another nearby BR system with relative error 
$4.98\%$ (using the weighted norm, with EAO-FGM), 
and with the following relative errors 
\begin{equation*}
\frac{{\|A- \tilde A\|}_F}{{\|A\|}_F} = 4.94\%,~
\frac{{\|B-  \tilde B\|}_F}{{\|B\|}_F} = 5.58\%,~
\frac{{\|C- \tilde  C\|}_F}{{\|C\|}_F} = 6.33\%,~
\frac{{\|D- \tilde D\|}_F}{{\|D\|}_F} = 5.12\%.
\end{equation*}
%Hence playing with the weight vector $w$ allows one to balance the importance of each matrix defining the system. 

\reviseELA{
\begin{remark}[Choice of the weights] 
The weights can be chosen by a trial-and-error approach, depending on the end goal of the user. For example, if the goal is to have similar relative errors for each matrix $(A,B,C,D)$ as shown above, one can use initial weights proportional to $(\|A\|_F^{-2}, \|B\|_F^{-2}, \|C\|_F^{-2}, \|D\|_F^{-2})$, and then adapts them iteratively. The algorithm can be initialized with the solution obtained with the previous weights to start from a better solution and hence have faster convergence. 
\end{remark}
}

\subsection{RLC ladder circuit}  \label{sec:karim} 

Let us now consider a larger system, with $n=200$ and $m=1$. It is an RLC ladder circuit with 100 resistors, inductors and capacitors~\cite{morGugA03}. We would like to obtain the sPH form of this system. The PH data-driven modeling of this system in  the impedance passive case was studied in~\cite{Che22}. First, we transform the system to the scattering representation from the impedance representation by applying the external Cayley transform~\cite{CheGHM23}. Given an impedance passive system ($A_I,B_I,C_I,D_I$), a scattering representation ($A_S,B_S,C_S,D_S$) can be obtained based on the following transformation 
\begin{equation*}
\begin{bmatrix}
A_S&B_S\\C_S&D_S
\end{bmatrix}=\begin{bmatrix}
A_I-B_I(I_m+D_I)^{-1}C_I&\sqrt{2}B_I(I_m+D_I)^{-1}\\-\sqrt{2}(I_m+D_I)^{-1}C_I&-(I_m+D_I)^{-1}(D_I-I_m) 
\end{bmatrix} .
\end{equation*}

\paragraph{Checking BRness} 

\begin{figure}[ht!]
\begin{center}
 \includegraphics[width=0.4\textwidth]{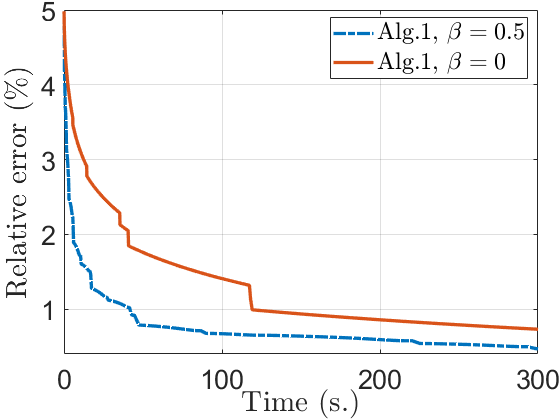}  
\caption{RLC ladder circuit with $n=200,m=1$: Relative error of Algorithm~\ref{alg:checkBRfgm} with extrapolation ($\beta = 0.5$) and without extrapolation ($\beta = 0$). 
%Note that the curve $(\beta=0.5,I=2)$ on the left corresponds to the curve $\rho = 10$ on the right. 
\label{fig:karimAlg1}}  
\end{center}
\end{figure}

Solving~\eqref{eq:mainoptprob2} with CVX provides a relative error of $1.3 \ 10^{-11}\%$ in 23.6 minutes, which shows that the system is BR, up to machine precision. This provides us with an sPH form. However, due to the numerically sensitive transformation $\hat Q = Q^{-1}$, this sPH form has a relative error of $6.73\ 10^{-4}\%$; see the next paragraph for an improved solution.  
Using a 5-minute timelimit, 
Algorithm~\ref{alg:checkBRfgm} finds a solution with relative error of  0.47\%. Figure~\ref{fig:karimAlg1} displays the evolution of the relative error over time, 
we compare the variant with extrapolation ($\beta = 0.5$) and without extrapolation ($\beta = 0$) to show the acceleration effect of extrapolation. 
Because Algorithm~\ref{alg:checkBRfgm} is a first-order method, it has a hard time computing a high-accuracy solution. If we run Algorithm~\ref{alg:checkBRfgm} for 25 minutes (for a fair comparison to CVX), we obtain a relative error of  0.13\%. This leads to an sPH form with relative error of 0.25\% (the increase is due to the change of variable $\hat Q = Q^{-1}$).

\paragraph{Nearest BR system in sPH form}  

\begin{figure}[ht!]
\begin{center}
 \includegraphics[width=0.5\textwidth]{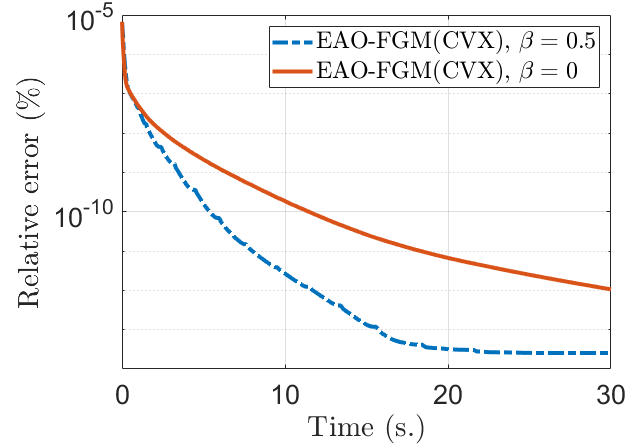}  
\caption{RLC ladder circuit with $n=200,m=1$: Relative error of Algorithm~\ref{alg:eao} with extrapolation ($\beta = 0.5$) and without extrapolation ($\beta = 0$), initialized by the CVX solution for problem~\eqref{eq:mainoptprob2}. 
%Note that the curve $(\beta=0.5,I=2)$ on the left corresponds to the curve $\rho = 10$ on the right. 
\label{fig:karimAlg2cvx}}  
\end{center}
\end{figure}

Using the solution obtained with CVX by solving~\eqref{eq:mainoptprob2}, we construct an sPH form with relative error $6.73\ 10^{-4}\%$. It can be further refined using Algorithm~\ref{alg:eao}. 
Figure~\ref{fig:karimAlg2cvx} displays the evolution of the relative error over time. Again, we compare the variants with extrapolation ($\beta = 0.5$) and without extrapolation ($\beta = 0$) to show the acceleration effect of extrapolation. Within 20 seconds, EAO-FGM produces a BR system in sPH form with relative error below $10^{-13} \%$.

\begin{figure}[ht!]
\begin{center}
 \includegraphics[width=0.5\textwidth]{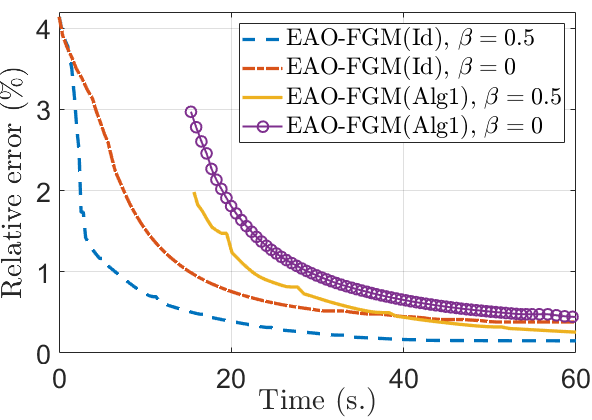}  
\caption{RLC ladder circuit with $n=200,m=1$: Relative error of several variants of Algorithm~\ref{alg:eao}, namely EAO-IPM(Id), with and without extrapolation, and EAO-FGM(Alg1) where Alg1 is given 15 seconds to initialize EAO-FGM.  
%Note that the curve $(\beta=0.5,I=2)$ on the left corresponds to the curve $\rho = 10$ on the right. 
\label{fig:karimAlg2}}  
\end{center}
\end{figure}

Figure~\ref{fig:karimAlg2} shows the evolution of the error of EAO-FGM(Id) and of EAO-FGM(Alg1), both with  and without extrapolation,  where Algorithm~\ref{alg:checkBRfgm} is given 25\% of the time to compute the initialization. We do not report the results for EAO-IPM which takes about 10 minutes to compute a single iteration. 
%using both CVX and FGM, and also Algorithm~\ref{alg:eao} initialized with the solution of Algorithm~\ref{alg:checkBRfgm} (from the previous section, time limit of 5 minutes). 
%As for the small system in the previous section, EAO-CVX is extremely slow (it performs one iteration in about 10 minutes): there is no reason to solve the subproblems in $Q$ and $(J,Z)$ exactly; this behavior will be confirmed for synthetic data in the next section. 
EAO-FGM(Id) with extrapolation performs best, and  produces a solution with a relative error of 0.15\% within one minute.   

\revise{
\begin{remark}[Initialization]
    \rm{
In the example of Figure~\ref{fig:karimAlg2}, the identity initialization performs better than Algorithm~\ref{alg:eao}. In fact, there is not guarantee that initializing with Algorithm~\ref{alg:eao} leads to better solutions, since we are dealing with hard non-convex problem. However, if the system is close to BRness, it is expected that  Algorithm~\ref{alg:eao} performs better since it would be able to find a solution with error close to zero. 
    }
\end{remark}
}

\subsection{Synthetic data} 

To provide a better idea of the performance of the various proposed algorithms, we generate synthetic data as follows: 
\begin{itemize}
    \item $F$, $P$ and $\tilde D$ are randomly generated using the normal distribution, and then normalized to have Frobenius norm 1. 
    In MATLAB, for the matrix $F$:  
    \texttt{F = randn(n,m); F = F/norm(F,'fro')}.
    Note that, in our code, we use $\tilde D = S+N$ where $S = \sym(\tilde{D})$ and $N = \skewp(\tilde D)$; this was the notation used for impedance passive systems in~\cite{gillis2018finding}, allowing us to use the same variables and a consistent PH representation. 

    \item The matrices $J$ is generated in the 
    same way, except that it is made skew-symmetric  before normalization. In MATLAB: 
    \texttt{J = randn(n); J = (J-J')/2; J = J/norm(J,'fro')}. 

\item For $Q$ and $R$, we first generate a matrix with randomly generated entries using the normal distribution, make the product with its transpose to obtain a PSD matrix, and then normalize it. 
In MATLAB: 
    \texttt{Q = randn(n); Q = Q*Q'; Q = Q/norm(Q,'fro')}. 

\item Finally, we construct $Z$ as in~\eqref{eq:defZ}. However, $Z$ might not be PSD. 
For any $c \geq -\lambda_{\min}(Z)$, where $\lambda_{\min}(Z)$ denotes the smallest eigenvalue of $Z$, $(Z+c I)/(1+c)$ is PSD while the (2,2) and (3,3) blocks remain equal to the identity matrix. 
We let $c = \max\left(0.1-\lambda_{\min}(Z),0.1\right)$, $R = (R + c \, I)/c$ while we divide by $1+c$ the off-diagonal entries of $Z$ (that is, $F$, $P$, $\tilde D$).  %{\color{red}{(PS: Replace $S$ and $N$ by })}
\end{itemize}

This procedure generates a synthetic BR system in sPH form. To make it non-BR, we add random noise to each matrix in the sPH form. More precisely, we generated Gaussian matrices, normalized them, possibly make them (skew-)symmetric when the corresponding factors are, 
and added them to each matrix multiplied by $\delta$. Since the generated matrices have norm 1, $\delta$ can be interpreted as a noise level.

\paragraph{Scalability of CVX to check BRness} 

Table~\ref{tab:timeCVXcheckBR} reports the average time and standard deviation for CVX to check BRness, that is, solve~\eqref{eq:mainoptprob2}, for various values of $n$ and $m$. The average is made on 10 randomly generated BR systems, and 10 randomly generated non-BR systems with noise level $\delta = 1$. 
\begin{center}
\begin{table}[h!]
\begin{center}
\label{numreal1}  
\begin{tabular}{c|c|c|c|c|c} 
 & $n= 10$, $m=5$ 
 & $n= 25$, $m=15$ 
 & $n= 50$, $m=25$
 & $n=75$, $m=40$  
 & $n=100$, $m=50$ \\  \hline
 noiseless & 
 $0.33 \pm 0.17$ s
 &  $1.28 \pm 0.35$ s 
 &  $29.9 \pm  6.2$ \; s  
 & $141.9 \pm  26.9$ s 
 & $8.3 \pm 1.1$ min \\
  noisy & 
 $0.30 \pm 0.03$ s
 &  $1.52 \pm 0.33$ s 
 &  $36.7 \pm  12.4$ s
 & $180.3  \pm 11.4$ s 
 & $11.1 \pm  0.35$ min
\end{tabular} 
\caption{Average time and standard deviation for SDPT3 (called via CVX) to solve~\eqref{eq:mainoptprob2} for 10 noiseless and 10 noisy ($\delta = 1$) data sets, except for $n=100$ for which we used only 5 datasets of each type, given the time needed to compute the solutions. \label{tab:timeCVXcheckBR}}
\end{center}
\end{table}
\end{center}

On a standard laptop, we can solve problems with size up to a hundred within a few minutes. For larger problems, we would need to resort to other strategies, such as Algorithm~\ref{alg:eao}.

\paragraph{How good is Algorithm~\ref{alg:eao} on synthetic data?}  

On the synthetic data without noise, that is, the generated system is BR, it takes about 8 minutes for CVX to solve~\eqref{eq:mainoptprob2} for $n=100$ and $m=50$; see Table~\ref{tab:timeCVXcheckBR}. 
Since the system is BR, this allows one to obtain a solution 
to~\eqref{eq:mainoptprob}. As we have seen in Section~\ref{sec:karim}, it is hard for Algorithm~\ref{alg:eao} to compute high-accuracy solutions in such cases. However, we may wonder how long it takes Algorithm~\ref{alg:eao}, on average, to compute a solution with a prescribed relative error. 

Table~\ref{tab:Alg2accur} reports the time needed for Algorithm~\ref{alg:eao} with the variant EAO-FGM(Id) to compute solutions of various accuracies as a function of time, averaged over 10 synthetic BR systems with $n=100$ and $m=50$. 
For example, in about 40 seconds on average, EAO-FGM(Id)   reaches solutions with relative error below 0.01\%, and, in less than 4 minutes on average, it reaches solutions with relative error below 0.001\%.  
\begin{center}
\begin{table}[h!]
\begin{center}
\begin{tabular}{c|c|c|c|c} 
relative error &    1\% 
 & 0.1\% 
 & 0.01\% 
 & 0.001\% 
 \\  \hline
time (s) &   
$0.53 \pm  0.02$ 
 & $6.2 \pm 0.3 $
 & $43.0 \pm  4.1$
 & $231 \pm  39$
 \\  \hline
 num. iter. &  
 $4 \pm   0$
 & $37 \pm 1.6$
 & $172   \pm   9.8 $
 & $369 \pm 18.8$
\end{tabular} 
\caption{Average time and number of iterations, and their standard deviation, to get a solution with a prescribed relative error for synthetic BR systems with $n=100$ and $m=50$. \label{tab:Alg2accur}}
\end{center}
\end{table}
\end{center}

 Table~\ref{tab:Alg2accurn200} reports the results for the same experiment but for $n=200$ and $m=100$. In this case, about 4 minutes are necessary to achieve a relative error below 0.01\%. Interestingly, the number of iterations to reach 0.01\% relative error is about the same as for the synthetic BR systems with $n=100$ and $m=50$. 
\begin{center}
\begin{table}[h!]
\begin{center}
\begin{tabular}{c|c|c|c} 
relative error &    1\% 
 & 0.1\% 
 & 0.01\% 
 \\  \hline
time (s) &   
$1.95 \pm  0.93$ 
 & $27.1 \pm 11.8 $
 & $253 \pm 99$
 \\  \hline
 num. iter. &  
 $3 \pm   0$
 & $27.6 \pm 1.3$
 & $160.1   \pm   3.7 $
\end{tabular} 
\caption{Average time and number of iterations, and their standard deviation, to get a solution with a prescribed relative error for synthetic BR systems with $n=200$ and $m=100$. \label{tab:Alg2accurn200}}
\end{center}
\end{table}
\end{center}

\paragraph{EAO-FGM vs.\ EAO-IPM}  

\begin{figure}[ht!]
\begin{center}
 \includegraphics[width=0.5\textwidth]{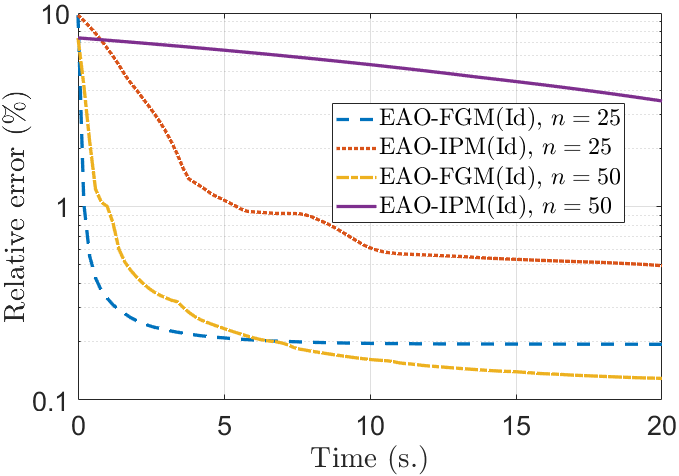}  
\caption{EAO-FGM vs.\ EAO-IPM: average relative error over time on 10 synthetic data sets with $(n=25,m=10)$ and $(n=50,m=25)$. 
\label{fig:FGMvsIPM}}  
\end{center}
\end{figure}

Let us compare again EAO-IPM and EAO-FGM, on synthetic data this time. 
Figure~\ref{fig:FGMvsIPM} shows the evolution over time of the average error for 10 randomly generated systems of size $(n=25,m=10)$ and $(n=50, m=25)$, with noise level $\delta = 0.1$.  
This confirms our observations from the previous sections: EAO-IPM spends a lot of time solving the subproblems exactly and is therefore significantly slower than EAO-FGM.

% \paragraph{Initialization}  

% As a rule of thumb, when the input system is close to being BR, then using CVX or Algorithm 1 to initialize it is highly beneficial. If the input system is far from being BR, then the identity initialization, which is very cheap, typically works better. 

% To illustrate this, Figure~\ref{} 
% shows the evolution the error ... average over 10 
%  $n=50$, $m=25$. We use a timelimit of 40 seconds, and run Algorithm~\ref{alg:checkBRfgm} for 10 seconds for the initialization. 
 
%  Id vs.\ FGM 

% \paragraph{Parameters} 

% $\rho$, $I$, 

% It seems rather sensitive to the input data... 

\section{Conclusion} 

In this paper, we studied the problem of finding the nearest bounded-real PH system to a given system. 
We derived the necessary conditions for a PH representation to be bounded real (Theorem~\ref{thm:suff_BR}). 
These conditions were used to derive a formulation to find the nearest scattering-passive port-Hamiltonian representation to a given system; see~\eqref{eq:mainoptprob}. We devised an algorithm that combines alternating optimization and Nesterov's fast gradient method (Algorithm~\ref{alg:eao}). In addition, we presented an algorithm to check if a system is bounded real using the same strategy (Algorithm~\ref{alg:checkBRfgm}). We showed the effectiveness of these strategies on real and synthetic systems. 
Further work includes learning bounded-real PH systems directly from data, and considering other error measures for the objective function, e.g., based on the transfer function. 

%Minimize the H2 norm of the transfer function? 

 %\section*{Statements}

\paragraph{Funding}  NG acknowledges the support by the European Union (ERC consolidator, eLinoR, no 101085607).  PS acknowledges the support of the SERB-CRG grant (CRG/2023/003221) and SERB-MATRICS grant by Government of India.

\reviseELA{
\section*{Acknowledgments}  We are grateful to the anonymous referee who
carefully read the manuscript; their feedback helped us improve our paper.
}

\small

\bibliographystyle{spmpsci}
\bibliography{Article2023}

\end{document}